\documentclass[3p,11pt,fleqn,numbers,square,sort&compress]{elsarticle}
\usepackage{amssymb,amsbsy,amsmath,amsfonts,amssymb,amscd}
\usepackage[english]{babel}
\usepackage{color}
\usepackage{float}
\usepackage{amssymb,natbib}
\usepackage{graphicx}
\usepackage[latin1]{inputenc}
\usepackage{times}
\usepackage{stmaryrd}
\usepackage{mathrsfs}
\usepackage{booktabs}
\usepackage{multirow,booktabs}
\usepackage{array}
\usepackage{subfig}
\usepackage{geometry}
\usepackage[ruled]{algorithm2e}
\usepackage[hidelinks,colorlinks=true,citecolor=blue,linkcolor=red,urlcolor=green,pdfstartview=]{hyperref}
\usepackage{pifont}
\usepackage{lineno}

\biboptions{comma,round}
\newtheorem{e-proposition}[theorem]{Proposition}

\newtheorem{e-definition}[theorem]{Definition\rm}

\def\og{\leavevmode\raise.3ex\hbox{$\scriptscriptstyle\langle\!\langle$~}}
\def\fg{\leavevmode\raise.3ex\hbox{~$\!\scriptscriptstyle\,\rangle\!\rangle$}}

\journal{Computers and Geotechnics}

\graphicspath{{fig/}}
\begin{document}

\begin{frontmatter}



\title{Numerical simulation of forerunning fracture in saturated porous solids with hybrid FEM/Peridynamic model}

\author[label1,label2]{Tao Ni}
\author[label3]{Francesco Pesavento}
\author[label1,label2]{Mirco Zaccariotto}
\author[label1,label2]{Ugo Galvanetto}
\author[label3,label4]{Bernhard A. Schrefler\corref{cor1}}

\ead{bernhard.schrefler@dicea.unipd.it}
\cortext[cor1]{Corresponding author}
%\address[label1]{School of Mechanics and Civil Engineering, China University of Mining and Technology, 221116, Xuzhou, China}
\address[label1]{Industrial Engineering Department, University of Padova, via Venezia 1, Padova, 35131, Italy}
\address[label2]{Center of Studies and Activities for Space (CISAS)-G. Colombo, via Venezia 15, Padova, 35131, Italy}
\address[label3]{Department of Civil, Environmental and Architectural Engineering, University of Padova, via Marzolo 9, Padova, 35131, Italy}
\address[label4]{Institute for Advanced Study, Technische Universität München, Lichtenbergstrasse 2a, D-85748 Garching b. München, Germany}

\begin{abstract}
In this paper, a novel hybrid FEM and Peridynamic modeling approach proposed in \citep{ni2020hybrid} is used to predict the dynamic solution of hydro-mechanical coupled problems.  A modified staggered solution algorithm is adopted to solve the coupled system. A one-dimensional dynamic consolidation problem is solved first to validate the hybrid modeling approach, and both $\delta-$convergence and $m_{r}-$convergence studies are carried out to determine appropriate discretization parameters for the hybrid model. Thereafter, dynamic fracturing in a rectangular dry/fully saturated structure with a central initial crack is simulated both under mechanical loading and fluid-driven conditions. In the mechanical loading fracture case, fixed surface pressure is applied on the upper and lower surfaces of the initial crack near the central position to force its opening. In the fluid-driven fracture case, the fluid injection is operated at the centre of the initial crack with a fixed rate. Under the action of the applied external force and fluid injection, forerunning fracture behavior is observed both in the dry and saturated conditions.
\end{abstract}


\end{frontmatter}

\section{Introduction}

The dynamic fracture propagation in saturated porous solids is not always
smooth and continuous %
\citep{cao2018porous,peruzzo2019dynamics,peruzzo2019stepwise}. The stepwise
fracture tip advancement and pressure oscillation in saturated porous media
have duly been documented with experiments in %
\citep{zhang2010dynamic,pizzocolo2013mode,razavi2016optimization,deng2016investigation}
and field observations in %
\citep{morita1990theory,fuh1992new,okland2002importance,fisher2012hydraulic,soliman2014analysis,de2015hydraulic}%
. In geophysical observations, intermittent fracture advancement is also
detected in saturated formations. Recognizing the existence of such a
behavior makes it easier to explain the non-volcanic (subduction) tremor and
volcanic tremor \citep{burlini2008volcanic,schwartz2007slow}. Furthermore,
properly predicting this phenomenon has also economic significance in
guiding the exploitation of resources in the reservoirs by using fracking
operations \citep{peruzzo2019stepwise,soliman2014analysis}.

A variety of numerical tools have been developed to simulate the dynamic
hydraulic fracturing in saturated media. A fully coupled cohesive-fracture
discrete model was combined to a generalized finite element formulation in %
\citep{schrefler2006adaptive} for solving the dynamic fracture problems in
porous media in thermal-hydro-mechanical coupled context. Further, the
three-dimensional version of the model was presented in %
\citep{secchi2012method} and applied to deal with the problem of hydraulic
fracturing in a concrete dam in quasi-static conditions. In %
\citep{jahromi2013development}, a three-dimensional, three-phase coupled
numerical model was introduced for hydro-mechanical coupled problems, taking
the mutual influence between dynamic fracture propagation and reservoir flow
into consideration. The double porosity approach was applied in %
\citep{kim2014fracture} to model the dynamic failure resulting from tensile
and shear stresses, dynamic non-linear permeability, leak-off and
thermal-poro-mechanical effects. In the numerical implementation of that
model, finite volume method was used for fluid flow, finite element method
for the solid phase and the backward Euler method for time discretization.
In \citep{ahn2014modeling}, the interactions between hydraulic fractures and
pre-existing natural fractures were investigated by using a coupled
numerical model that integrates dynamic fracture propagation and reservoir
flow simulation. A hydraulic fracture simulator for non-isothermal
conditions was developed in \citep{kim2013development} by coupling a flow
simulator to a geomechanics code. The developed simulator was then applied
to dynamic hydraulic fracture cases in \citep{kim2015numerical}, in which
the fracture advancement was found to occur discontinuously in time and the
pressure showed a saw-toothed response. A coupled flow and geomechanics
model was developed in \citep{feng2017parameters}, where a cohesive zone
model was adopted for simulating fluid-driven fracture and a
poro-elastic/plastic formation behavior. In the above-mentioned works, all
solutions featured irregular fracture advancement steps which point to a
physical origin, generally documented by the fact that there are fewer crack
advancement steps in a time interval than number of time steps. However, in
some other works %
\citep{mohammadnejad2013hydro,mohammadnejad2013extended,fathima2019implications}%
, the fracture showed a regular advancement process with the adopted
simulation parameters. As discussed in \citep{peruzzo2019stepwise}, regular
fracture advancement steps are mainly of numerical origin most probably
caused by improper time step/fracture advancement algorithms %
\citep{secchi2014hydraulic,cao2017interaction}.

Besides the observed fracture propagation with irregular steps, forerunning
fracturing is a further proof of the existence of stepwise advancement %
\citep{peruzzo2019stepwise}. Forerunning fracture behavior exists in rocks %
\citep{sammonds1989acoustic} and has been observed in the numerical
simulation of a double cantilever beam problem in \citep{peruzzo2019stepwise}%
, under the action of a pre-defined monotonically increasing load.
Forerunning was also observed in \citep{slepyan2015forerunning} and in case
of heterogeneous material in \citep{tvergaard1993analysis}. These examples
took only mechanical wave propagation into consideration, ignoring the
effects of the wave propagation of the pore pressure in a saturated porous
medium. No other work on the numerical simulation of forerunning fracture in
saturated porous media has been found in the literature. This motivates
further investigation on this issue in this paper, where a hybrid
FEM/Peridynamic approach will be applied.

Peridynamics (PD), firstly introduce by Silling in 2000 %
\citep{silling2000reformulation}, is a non-local continuum theory based on
spatial integro-differential equations, which has an unparalleled capability
to simulate the crack propagation in structures due to allowing cracks to
grow naturally without resorting to external crack growth criteria. The
first version of the PD theory, called Bond-Based PD (BB-PD), had a strong
limitation because the Poisson's ratio could only assume a fixed value %
\citep{silling2005meshfree,zhu2017peridynamic,wang2018three}. To remove
that limitation, its most recent version named as state-based PD (SB-PD) is
introduced in \citep{silling2007peridynamic}, including the ordinary and
non-ordinary versions (OSB-PD and NOSB-PD). In the past two decades,
PD-based computational methods have made great progress in dealing with
crack propagation problems %
\citep{silling2005meshfree,lai2015peridynamics,bobaru2016handbook,cheng2019dynamic,zhang2019failure,diana2020simulating}%
. Although the PD-based numerical approaches have great advantages in
solving crack propagation problems, they are usually more computationally
expensive than those making use of local mechanics and FEM %
\citep{galvanetto2016effective,zaccariotto2018coupling,ni2019static}.

Inspired by \citep{milanese2016avalanches} and with the aim to reduce the
computing costs, a novel hybrid FEM and Peridynamic modeling approach has
been proposed in \citep{ni2020hybrid} to simulate the hydro-mechanical
coupled fracturing process in saturated porous media. This hybrid modeling
approach is here used to solve dynamic hydro-mechanical coupled problems.
The staggered solution algorithm is adopted for this purpose. In each
solving sequence, the implicit time integration iteration from %
\citep{zienkiewicz2000finite} is used to solve the FE equations of fluid
flow, and a modified explicit central difference time integration algorithm
presented in \citep{taylor1989pronto} is used to solve the peridynamic
equations. A one-dimensional dynamic consolidation problem is first
addressed for validation of the presented approach. Subsequently, the
dynamic fracture propagation in a rectangular porous structure with a
central initial crack under mechanical loading and under fluid injection is
simulated to investigate the possibility of forerunning fracture behavior.
As mentioned in \citep{slepyan2015forerunning}, the action of high-amplitude
incident sinusoidal waves could lead to forerunning mode of fracture.
Therefore, the external force and fluid injection rate applied to the
rectangle porous samples are selected with significantly large values such
as to generate a high-amplitude wave in the studied system.

\section{Methodology}

\subsection{Ordinary state-based peridynamic formulation for fully saturated
porous media}

\begin{figure}[h]
\begin{center}
\includegraphics[scale=0.5]{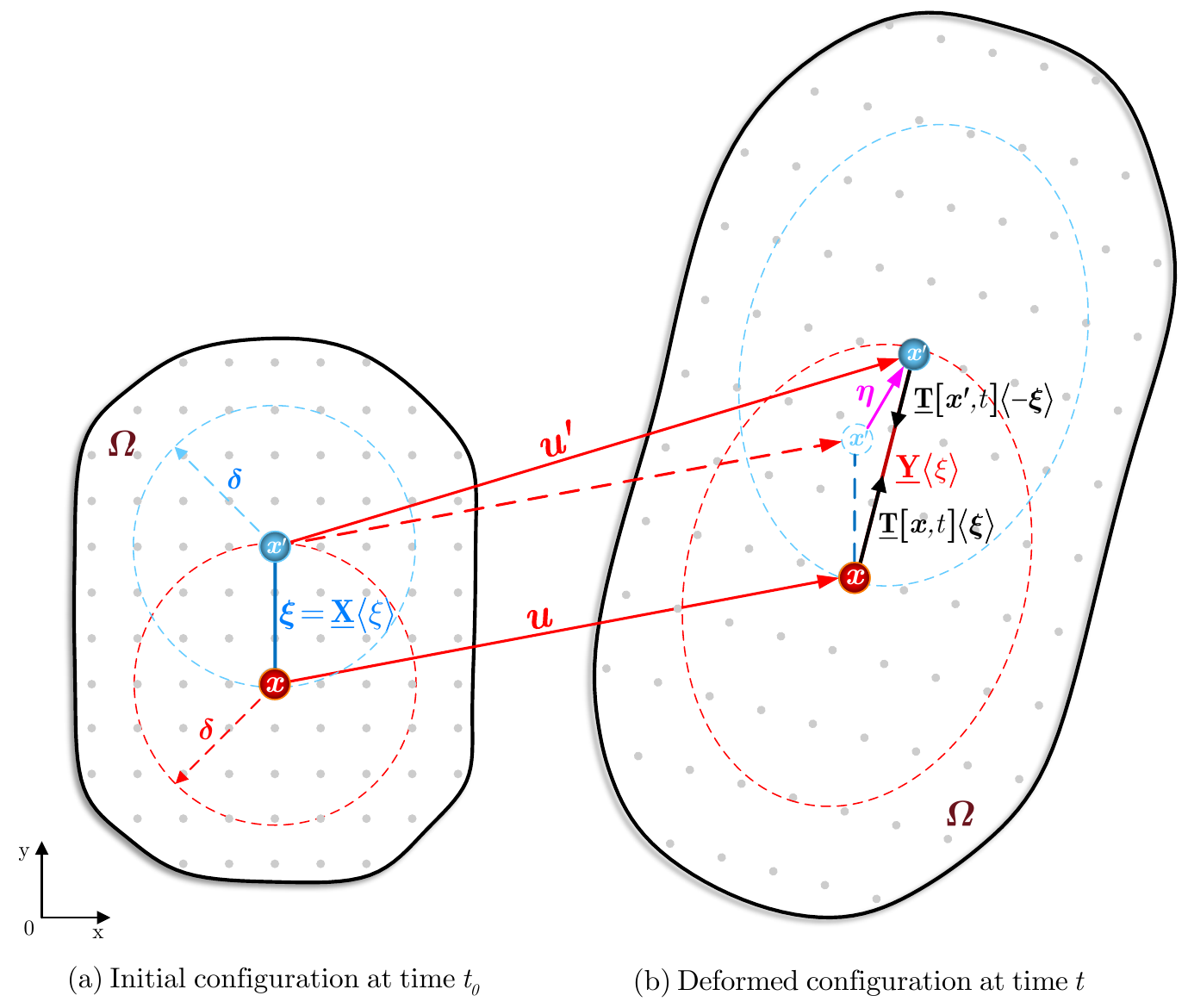}
\end{center}
\caption{The positions of two material points in the (a) initial and (b)
deformed configurations.}
\label{fig1}
\end{figure}
As shown in Fig. \ref{fig1}, a domain $\Omega $ modeled by the OSB-PD is
composed by material points of infinitesimal volume, where each point would
interact with the other points around it within a prescribed horizon radius (%
$\delta $). Assuming that there is a point $\boldsymbol{x}$ in a domain $%
\Omega $, and $\boldsymbol{x}^{\prime }$ is a point within its horizon, the
hydro-mechanical coupled peridynamic equation of motion at time $t$ can be
given as \citep{ni2020hybrid}: 
\begin{equation}
\rho \boldsymbol{\ddot{u}}\left( \boldsymbol{x},t\right) =\int\nolimits_{ 
\mathcal{H}_{x}}\left\{ \text{\textbf{\b{T}}}\left[ \boldsymbol{x},p,t\right]
\left\langle \boldsymbol{\xi }\right\rangle -\text{ \textbf{\b{T}}}\left[ 
\boldsymbol{x}^{\prime },p^{\prime },t\right] \left\langle \boldsymbol{-\xi }%
\right\rangle \right\} dV_{x^{\prime }}+\boldsymbol{b}\left( \boldsymbol{x}%
,t\right)  \label{2.1}
\end{equation}%
where $\rho$ is the material mass density, $\boldsymbol{\ddot{u}}$ is the
acceleration, $\text{\textbf{\b{T}}}\left[ \boldsymbol{x},p,t\right]
\left\langle \boldsymbol{\xi }\right\rangle $ and $\text{ \textbf{\b{T}}}%
\left[ \boldsymbol{x}^{\prime },p^{\prime },t\right] \left\langle 
\boldsymbol{-\xi }\right\rangle$ are force density vector states along the
deformed bond at points $\boldsymbol{x}$ and $\boldsymbol{x}^{\prime }$,
respectively. $dV_{x^{\prime }}$ is the infinitesimal volume associated to
point $\boldsymbol{x}^{\prime }$, $\boldsymbol{b}$ is the force density
applied by some external force. To describe the position of $\boldsymbol{x}%
^{\prime }$ with respect to $\boldsymbol{x}$, a relative position vector is
defined as: 
\begin{equation}
\boldsymbol{\xi =x}^{\prime}-\boldsymbol{x}  \label{2.2}
\end{equation}

After being displaced respectively by $\boldsymbol{u}$ and $\boldsymbol{u}%
^{\prime }$, the relative displacement vector of $\boldsymbol{x}$ and $%
\boldsymbol{x}^{\prime }$ is given as: 
\begin{equation}
\boldsymbol{\eta }=\boldsymbol{u}^{\prime }-\boldsymbol{u}  \label{2.3}
\end{equation}

To specifically describe the initial and deformed states of the bond in the
OSB-PD theory, the reference position vector state, $\underline{\boldsymbol{X%
}}\left\langle \boldsymbol{\xi }\right\rangle =\boldsymbol{\xi }$, and the
deformation vector state, $\underline{\boldsymbol{Y}}\left\langle 
\boldsymbol{\xi }\right\rangle =\boldsymbol{\xi +\eta }$, are defined.
Therefore, the reference position scalar state and the deformation scalar
state, representing the lengths of the bond in its initial and deformed
state, respectively, are given as: 
\begin{equation}
\begin{array}{ccc}
\underline{x}=\left\Vert \underline{\boldsymbol{X}}\right\Vert & , & 
\underline{y}=\left\Vert \underline{\boldsymbol{Y}}\right\Vert%
\end{array}
\label{2.4}
\end{equation}
where $\left\Vert \cdot \right\Vert $ denotes the Euclidean norm.

In a linear OSB-PD material, the elastic strain energy density at point $%
\boldsymbol{x}$ is given as \citep{silling2010linearized,bobaru2016handbook}: 
\begin{equation}
W\left( \theta ,\underline{e}^{d}\right) =\frac{k^{\prime }\theta ^{2}}{2}+%
\frac{\gamma }{2}\int\nolimits_{\mathcal{H}_{x}}\underline{\mathit{w}}\text{ 
}\underline{e}^{d}\text{ }\underline{e}^{d}\text{d}V_{x^{\prime }}
\label{2.5}
\end{equation}%
where $\theta $\ is the volume dilatation value, $\underline{e}^{d}$ is the
deviatoric extension state, $k^{\prime }$ and $\gamma $ are positive
constants related to mechanical material parameters, $\underline{\mathit{w}}$
is an influence function and several examples of its forms have been
summarized in \citep{ni2019coupling}.

For simplification, only plane strain condition is considered in this paper. 
$\theta $ and $\underline{e}^{d}$ are usually given as \citep{ni2019coupling}%
: 
\begin{equation}
\theta =\frac{2}{m}\int\nolimits_{\mathcal{H}_{x}}\left( \underline{\mathit{w%
}}\text{ }\underline{x}\text{ }\underline{e}\right) \text{d}V_{x^{\prime }}
\label{2.6}
\end{equation}%
\begin{equation}
\underline{e}^{d}=\underline{e}-\frac{\theta \underline{x}}{3}  \label{2.7}
\end{equation}%
in which, $\underline{e}$ is the extension scalar state defined as: 
\begin{equation}
\underline{e}=\underline{y}-\underline{x}  \label{2.8}
\end{equation}%
In plane strain cases, $k^{\prime }$ and $\gamma $ are given as: 
\begin{equation}
k^{\prime }=\kappa +\frac{\mu }{9}  \label{2.9}
\end{equation}%
\begin{equation}
\gamma =\frac{8\mu }{m}  \label{2.10}
\end{equation}%
where $\kappa $ and $\mu $ are bulk modulus and
shear modulus, respectively. $m$ is the weighted volume defined as: 
\begin{equation}
m=\int\nolimits_{\mathcal{H}_{x}}\underline{\mathit{w}}\left\Vert 
\boldsymbol{\xi }\right\Vert ^{2}\text{d}V_{x^{\prime }}  \label{2.11}
\end{equation}

\textcolor{blue}{In order to consider the effect of pore pressure on solid deformation, the effective stress principle \citep{lewis1998finite} need to be introduced in the peridynamic model. In \citep{song2019peridynamics,song2020peridynamic}, an effective force state is proposed in the framework of non-ordinary state-based peridynamics (NOSB-PD) to model geomaterials. The NOSB-PD model is indeed a good choice to simulate the complex behavior of solid materials by introducing the developed classical continuous constitutive relations. In this paper, the ordinary state-based peridynamic model is adopted, and according to \citep{turner2013non,ni2020hybrid},} the effective force density vector
state $\text{\textbf{\b{T}}}\left[ \boldsymbol{x},p,t\right] \left\langle 
\boldsymbol{\xi }\right\rangle$ in Eq.(\ref{2.1}) is defined as: 
\begin{equation}
\text{\textbf{\b{T}}}\left[ \boldsymbol{x},p,t\right] \left\langle 
\boldsymbol{\xi }\right\rangle=\underline{t}\left[ \boldsymbol{x},p,t\right]
\left\langle \boldsymbol{\xi }\right\rangle \cdot \underline{\boldsymbol{M}}%
\left\langle \boldsymbol{\xi }\right\rangle  \label{2.12}
\end{equation}
where $\underline{t}$ is called the force density scalar state. $\underline{%
\boldsymbol{M}}$ is a unit state in the direction of the deformed bond,
which is usually defined as: 
\begin{equation}
\underline{\boldsymbol{M}}\left\langle \boldsymbol{\xi }\right\rangle =\frac{
\underline{\boldsymbol{Y}}}{\left\Vert \underline{\boldsymbol{Y}}\right\Vert 
}  \label{2.13}
\end{equation}

Referring to \citep{turner2013non,ni2020hybrid}, $\underline{t}\left[ 
\boldsymbol{x},p,t\right] $ is defined in plane strain cases as: 
\begin{equation}
\underline{t}\left[ \boldsymbol{x},p,t\right] =\left[ \left( 2k^{\prime }-%
\frac{1}{9}\gamma m\right) \theta -2\alpha p\right] \frac{\underline{\mathit{%
w}}\text{ }\underline{x}}{m}+\gamma \underline{\mathit{w}}\text{ }\underline{%
e}^{d}=\underline{t}\left[ \boldsymbol{x},t\right] -2\alpha p\frac{%
\underline{\mathit{w}}\text{ }\underline{x}}{m}  \label{2.14}
\end{equation}

Thus, the hydro-mechanical coupled peridynamic equation of motion will become:%
\begin{equation}
\begin{array}{l}
\rho \boldsymbol{\ddot{u}}\left( \boldsymbol{x},t\right) =\int\nolimits_{ 
\mathcal{H}_{x}}\left\{ \text{\textbf{\b{T}}}\left[ \boldsymbol{x},p,t\right]
\left\langle \boldsymbol{\xi }\right\rangle -\text{ \textbf{\b{T}}}\left[ 
\boldsymbol{x}^{\prime },p^{\prime },t\right] \left\langle \boldsymbol{-\xi }%
\right\rangle \right\} dV_{x^{\prime }}+\boldsymbol{b}\left( \boldsymbol{x}%
,t\right) \\ 
=\int\nolimits_{\mathcal{H}_{x}}\left\{ \text{\textbf{\b{T}}}\left[ 
\boldsymbol{x},t\right] \left\langle \boldsymbol{\xi }\right\rangle -\text{ 
\textbf{\b{T}}}\left[ \boldsymbol{x}^{\prime },t\right] \left\langle 
\boldsymbol{-\xi }\right\rangle \right\} dV_{x^{\prime }}-2\alpha
\int\nolimits_{\mathcal{H}_{x}}\left[ p\frac{\underline{\mathit{w}}\text{ } 
\underline{x}}{m\left( \boldsymbol{x}\right) }\underline{\boldsymbol{M}}%
\left\langle \boldsymbol{\xi }\right\rangle -\text{ }p^{\prime }\frac{ 
\underline{\mathit{w}}\text{ }\underline{x}}{m\left( \boldsymbol{x}^{\prime
}\right) }\underline{\boldsymbol{M}}\left\langle -\boldsymbol{\xi }%
\right\rangle \right] dV_{x^{\prime }}+\boldsymbol{b}\left( \boldsymbol{x}%
,t\right)%
\end{array}
\label{2.14_1}
\end{equation}
where $p$ and $p^{\prime }$\ are the values of pore pressure at points $%
\boldsymbol{x}$ and $\boldsymbol{x}^{\prime }$.

\subsection{Failure criteria}

Failure criteria are essential in PD-based numerical models to describe
material failure and crack advancement. In \citep{zhang2018state}, a strain
energy-based \textquotedblleft critical bond stretch\textquotedblright\
failure criterion is introduced for state-based PD model and successfully
applied to the quantitative fracture analysis of solid materials. Here we
adopt that failure criterion for convenience.

As in \citep{zhang2018state}, when the stretch value of a bond reaches the
critical value $s_{c}$, it will be broken, indicating the advancement of
fracture. The stretch value of bond $\boldsymbol{\xi }$ is defined as %
\citep{ni2019coupling}: 
\begin{equation}
s\left\langle \boldsymbol{\xi }\right\rangle =\frac{\underline{e}%
\left\langle \boldsymbol{\xi }\right\rangle }{\underline{x}\left\langle 
\boldsymbol{\xi }\right\rangle }  \label{2.15}
\end{equation}%
thus, the extension scalar state can be also expressed as: 
\begin{equation}
\underline{e}\left\langle \boldsymbol{\xi }\right\rangle =s\left\langle 
\boldsymbol{\xi }\right\rangle \text{ }\underline{x}\left\langle \boldsymbol{%
\xi }\right\rangle  \label{2.16}
\end{equation}

\begin{figure}[h]
\begin{center}
\includegraphics[scale=0.5]{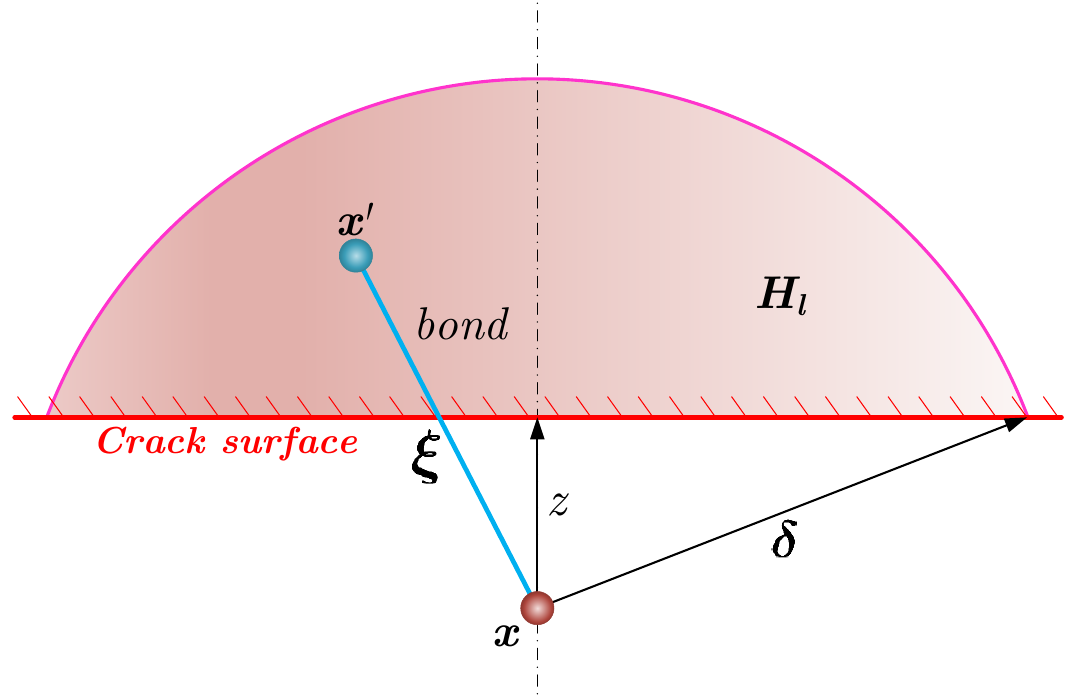}
\end{center}
\caption{A diagrammatic sketch of a peridynamic domain crossed by a crack
surface.}
\label{fig2}
\end{figure}

As shown in Fig. \ref{fig2}, $\mathcal{H}_{l}$ is the domain removed by the
crack surface from the neighborhood of point $\boldsymbol{x}$. $\boldsymbol{x%
}^{\prime }$ represents points located in the domain $\mathcal{H}_{l}$. The
bond connecting $\boldsymbol{x}$ and $\boldsymbol{x}^{\prime }$ is marked as 
$\boldsymbol{\xi }$. At the formation of the crack surface, the deformation
energy stored in the bond $\boldsymbol{\xi }$ is released. Therefore, we
could assume that the work required to break all the bonds connecting point $%
\boldsymbol{x}$ to points in the domain $\mathcal{H}_{l}$ is equal to the
summation of the energy stored in these bonds in their critical stretch
condition. In critical condition, the extension scalar states of these bonds
are $\underline{e}_{c}=s_{c}$ $\underline{x}$. 

According to the definition of $\theta $ in Eq. (\ref{2.6}), the volume
dilatation value caused by the broken bonds can be computed by: 
\begin{equation}
\theta _{\boldsymbol{l}}=\frac{2}{m}\int\nolimits_{\mathcal{H}_{l}}\left( 
\underline{\mathit{w}}\text{ }\underline{x}\text{ }\underline{e}_{c}\right) 
\text{ d}V_{x^{\prime }}=\frac{2}{m}\int\nolimits_{\mathcal{H}_{l}}\left( 
\underline{\mathit{w}}\text{ }\underline{x}^{2}\text{ }s_{c}\right) \text{ d}%
V_{x^{\prime }}=2\varpi _{l}\text{ }s_{c}  \label{2.17}
\end{equation}%
where 
\begin{equation}
\begin{array}{ccc}
  \varpi _{l}=\frac{m_{l}}{m}& \& & m_{l}=\int\nolimits_{\mathcal{H}_{l}}\left( \underline{\mathit{w}}\text{ }%
  \underline{x}^{2}\right) dV_{x^{\prime }}
\end{array}
\label{2.18}
\end{equation}
is the volume weight from the points in the domain $\mathcal{H}_{l}$, and $%
\varpi _{l}=\frac{m_{l}}{m}$ is the ratio of $m_{l}$\ to $m$.

Using Eqs. (\ref{2.5}), (\ref{2.6}), (\ref{2.9}) and (\ref{2.10}), the
expression of strain energy density can be rewritten as %
\citep{zhang2018state,silling2010linearized}: 
\begin{equation}
W\left( \theta ,\underline{e}\right) =\frac{1}{2}\left( \kappa -\frac{7}{9}%
\mu \right) \theta ^{2}+\frac{4\mu }{m}\int\nolimits_{\mathcal{H}_{x}}%
\underline{\mathit{w}}\text{ }\underline{e}\text{ }\underline{e}\text{ }%
\text{d}V_{x^{\prime }}  \label{2.19}
\end{equation}

Accordingly, substitution of Eqs. (\ref{2.17}-\ref{2.18}) into (\ref{2.19})
gives the energy released by the broken bonds at point $\boldsymbol{x}$ as: 
\begin{equation}
W\left\langle \boldsymbol{x}\right\rangle =\frac{1}{2}\left( \kappa -\frac{7%
}{9}\mu \right) \theta _{\boldsymbol{l}}^{2}+\frac{4\mu }{m}\int\nolimits_{%
\mathcal{H}_{l}}\underline{\mathit{w}}\text{ }\underline{x}^{2}\text{ }%
s_{c}^{2}\text{d}V_{x^{\prime }}=\left[ 2\varpi _{l}^{2}\left( \kappa -\frac{%
7}{9}\mu \right) +4\mu \varpi _{l}\right] s_{c}^{2}\text{ }  \label{2.20}
\end{equation}

Consequently, the released energy per unit fracture area is: 
\begin{equation}
G_{c}=2\int\nolimits_{0}^{\delta }W\left\langle \boldsymbol{x}\right\rangle 
\text{ d}z=2\int\nolimits_{0}^{\delta }\left[ 2\varpi _{l}^{2}\left( \kappa -%
\frac{7}{9}\mu \right) +4\mu \varpi _{l}\right] s_{c}^{2}\text{ d}z=\left[
4\left( \kappa -\frac{7}{9}\mu \right) \beta +8\mu \beta ^{\prime }\right]
s_{c}^{2}  \label{2.21}
\end{equation}

In Eq. (\ref{2.21}), $\beta $ and $\beta ^{\prime }$ are parameters affected
by the influence function, and are defined as: 
\begin{equation}
\begin{array}{ccc}
\beta =\int\nolimits_{0}^{\delta }\varpi _{l}^{2}\text{ d}z & , & \beta
^{\prime }=\int\nolimits_{0}^{\delta }\varpi _{l}\text{ d}z%
\end{array}
\label{2.22}
\end{equation}

%The values of $\beta $\ and $\beta ^{\prime }$\ are defined by the influence
%function $\underline{\mathit{w}}$ and horizon radius $\delta $. 
In this paper, the influence function is taken as $\underline{\mathit{w}}=1$%
, therefore, the values of $\beta $\ and $\beta ^{\prime }$\ can be obtained
as: 
\begin{equation}
\begin{array}{ccc}
\beta =0.2192\delta & , & \beta ^{\prime }=\frac{2\delta }{5\pi }%
\end{array}
\label{2.23}
\end{equation}%
According to Eq. (\ref{2.21}), the critical stretch value $%
s_{c}$ of plane strain OSB-PD model can be given as:

\begin{equation}
s_{c}=\sqrt{\frac{G_{c}}{\left[ 4\left( \kappa -\frac{7}{9}\mu \right)
\beta +8\mu \beta ^{\prime }\right] }}  \label{2.24}
\end{equation}
where $G_{c}$ is the critical energy release rate for mode I fracture. The
effect of damage can be taken into account by introducing a variable $%
\mathit{\varrho }$, in a way similar to that used in BB-PD models %
\citep{zaccariotto2018coupling,ni2018peridynamic}: 
\begin{equation}
\underline{\mathit{\varrho }}\left\langle \boldsymbol{\xi }\right\rangle 
\boldsymbol{=}\left\{ 
\begin{array}{ccc}
1 & , & \text{if }s\left\langle \boldsymbol{\xi }\right\rangle <s_{c}\text{
, for all }0<\overline{t}<t \\ 
0 & , & \text{otherwise}%
\end{array}%
\right.  \label{2.25}
\end{equation}%
then the damage value $\varphi _{x}$ at point $\boldsymbol{x}$\ in the
system can be defined as: 
\begin{equation}
\varphi _{x}=1-\frac{\int\nolimits_{\mathcal{H}_{x}}\underline{\mathit{w}}
\left\langle \boldsymbol{\xi }\right\rangle \text{ }\underline{\mathit{\
\varrho }}\left\langle \boldsymbol{\xi }\right\rangle \text{d}V_{x^{\prime }}%
}{\int\nolimits_{\mathcal{H}_{x}}\underline{\mathit{w}}\left\langle 
\boldsymbol{\xi }\right\rangle \text{d}V_{x^{\prime }}}  \label{2.26}
\end{equation}%
where $\varphi _{x}\in \left[ 0,1\right] $, and the cracks can be identified
wherever $\varphi _{x}\geqslant 0.5$.

\subsection{Formulation for fluid flow in fractured porous media}

\begin{figure}[h!]
\centering  
\subfloat[]{\includegraphics[scale=0.4]{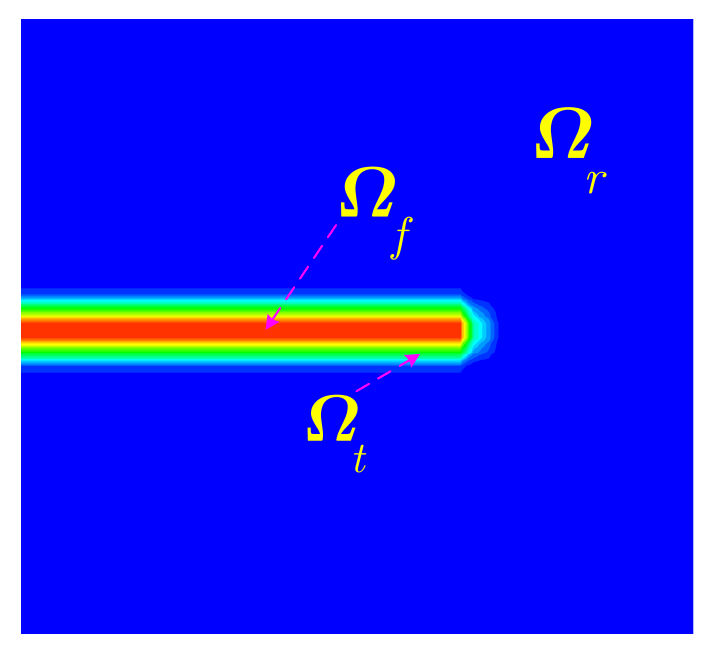}			%
\label{fig2_1:sub1}}\hspace{0.5in} \subfloat[]{		%
\includegraphics[scale=0.4]{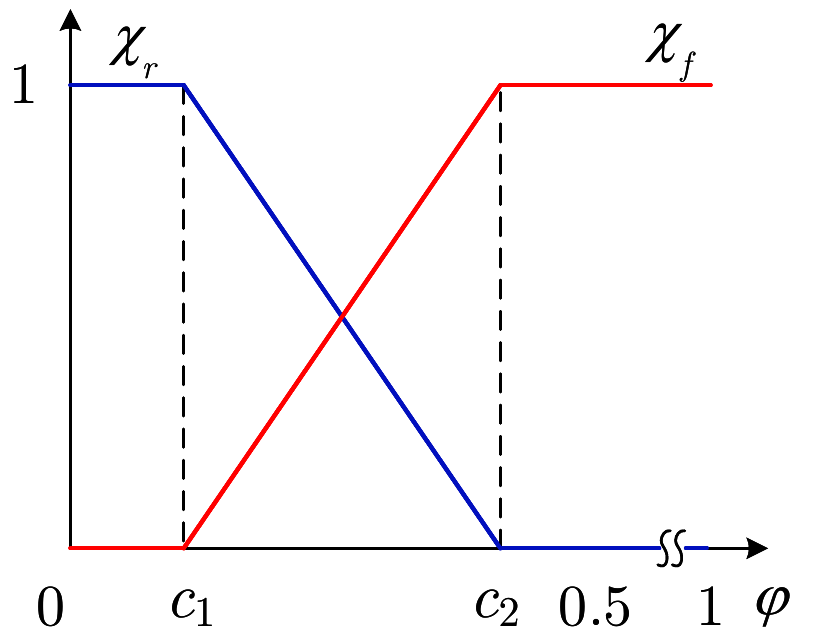}			%
\label{fig2_1:sub2}}
\caption{(a) The definition of the $\Omega _{r}$, $\Omega _{t}$ and $\Omega
_{f}$; (b) linear indicator functions $\protect\chi _{r}$ and $\protect\chi %
_{f}$.}
\label{fig2_1}
\end{figure}

In this paper, we use Darcy's law to describe the flow field in the
saturated porous domain $\Omega $. Other more complicated flow field models
will be considered in future work. To formulate the governing equation for
flow in the fractured porous domain, the whole domain $\Omega $ is divided
into three parts: $\Omega _{r}$, $\Omega _{f}$\ and $\Omega _{t}$ (see Fig. %
\ref{fig2_1:sub1}), representing the unbroken domain (reservoir domain), the
fracture domain and the transition domain between $\Omega _{r}$ and $\Omega
_{f}$. The peridynamic damage field ($\varphi $ in Eq. (\ref{2.26})) is used
as an indicator. As shown in Fig. \ref{fig2_1:sub2}, two threshold values ($%
c_{1}$ and $c_{2}$) are set to identify the three flow domains: the
reservoir domain is defined as $\varphi \leqslant c_{1}$, the fracture
domain as $\varphi \geqslant c_{2}$ and the transition domain as $c_{1}\leq
\varphi \leq c_{2}$.

Therefore, the mass balance equation for the reservoir domain is given as:%
\begin{equation}
\rho _{r}s_{r}\frac{\partial p}{\partial t}+\rho _{r}\alpha _{r}\frac{
\partial \varepsilon _{v}}{\partial t}+\rho _{r}\nabla \cdot \left[ \frac{
k_{r}}{\mu_{w}}\left( -\nabla p+\rho _{r}g\right) \right] =q_{r}
\label{2.2.1}
\end{equation}%
where $\alpha _{r}$,\ $s_{r}$, $k_{r}$, $q_{r}$ and $\rho _{r}$\ are the
Biot coefficient, storage coefficient, permeability, source term and the
density of the medium in the reservoir domain, respectively; $\mu_{w}$ is
the viscosity coefficient of the fluid in the reservoir domain. $g$\ is
gravity and $\varepsilon _{v}$\ is volumetric strain. The storage
coefficient is given as: 
\begin{equation}
s_{r}=\frac{\left( \alpha _{r}-n_{r}\right) \left( 1-\alpha _{r}\right) }{
K_{r}}+\frac{n_{r}}{K_{w}}  \label{2.2.2}
\end{equation}%
where $K_{r}$ and $K_{w}$ are the bulk moduli of solid and fluid in the
reservoir domain, $n_{r}$\ is the porosity.

The fracture domain is assumed to be fully filled by fluid and the porosity $%
n_{f}=1$ is adopted. Supposing that the fluid in the fracture is
incompressible, the volumetric strain term in Eq. (\ref{2.2.1}) vanishes %
\citep{van2019monolithic,lee2016pressure}, and the mass balance in the
fracture domain is expressed as: 
\begin{equation}
\rho _{f}s_{f}\frac{\partial p}{\partial t}+\rho _{f}\nabla \cdot \left[ 
\frac{k_{f}}{\mu_{w}}\left( -\nabla p+\rho _{f}g\right) \right] =q_{f}
\label{2.2.3}
\end{equation}%
where $s_{f}$, $k_{f}$, $q_{f}$ and $\rho _{f}$\ are the storage
coefficient, permeability, source term and the density of the fluid in the
fracture domain.

Two linear indicator functions are defined to connect the transition domain
with the reservoir and fracture domains: $\chi _{r}$ and $\chi _{f}$ %
\citep{lee2016pressure,zhou2018phase,zhou2019phase}, and they satisfy: 
\begin{equation}
\chi _{r}\left( \cdot ,\varphi \right) :=\chi _{r}\left( \boldsymbol{x}%
,t,\varphi \right) =1\text{ in }\Omega _{r}\left( t\right) \text{ and }\chi
_{r}\left( \cdot ,\varphi \right) =0\text{ in }\Omega _{f}\left( t\right)
\label{2.2.4}
\end{equation}
\begin{equation}
\chi _{f}\left( \cdot ,\varphi \right) :=\chi _{f}\left( \boldsymbol{x}%
,t,\varphi \right) =1\text{ in }\Omega _{f}\left( t\right) \text{ and }\chi
_{f}\left( \cdot ,\varphi \right) =0\text{ in }\Omega _{r}\left( t\right)
\label{2.2.5}
\end{equation}

As described in Fig. \ref{fig2_1:sub2}, the two linear functions are defined
as:%
\begin{equation}
\chi _{r}\left( \cdot ,\varphi \right) =\frac{c_{2}-\varphi }{c_{2}-c_{1}},%
\text{ and }\chi _{f}\left( \cdot ,\varphi \right) =\frac{\varphi -c_{1}}{
c_{2}-c_{1}}  \label{2.2.6}
\end{equation}

Fluid and solid properties in the transition domain can be obtained by
interpolating those properties of reservoir and fracture domains by using
the indicator functions: $\rho _{T}=\rho _{r}\chi _{r}+\rho _{f}\chi _{f}$, $%
\alpha _{T}=\alpha _{r}\chi _{r}+\alpha _{f}\chi _{f}$, $n_{T}=n_{r}\chi
_{r}+n_{f}\chi _{f}$, $k_{T}=k_{r}\chi _{r}+k_{f}\chi _{f}$, $%
s_{T}=s_{r}\chi _{r}+s_{f}\chi _{f}$. It should be noted in particular that
the Biot coefficient in the fracture domain is taken as $\alpha_{f}=1$.
Therefore, the governing equation for the fluid flow in the transition
domain can be given as:%
\begin{equation}
\rho _{T}s_{T}\frac{\partial p}{\partial t}+\rho _{T}\alpha _{T}\frac{
\partial \varepsilon _{v}}{\partial t}+\rho _{T}\nabla \cdot \left[ \frac{
k_{T}}{\mu_{w}}\left( -\nabla p+\rho _{T}g\right) \right] =q_{T}
\label{2.2.7}
\end{equation}

Moreover, the cubic law is used to evaluate the permeability in the fracture
domain \citep{zimmerman1994hydraulic,cao2017interaction}:%
\begin{equation}
k_{f}=\frac{1}{12}a^{2}  \label{2.2.8}
\end{equation}%
where $a$ is the aperture of the crack, and it can be obtained by using the
method described in the $sect.$ $3.3$ of \citep{ni2020hybrid}.

\section{Discretization and numerical solution procedure}

\begin{figure}[h!]
\centering  
\subfloat[Effect of fluid on
	solid.]{\includegraphics[scale=0.55]{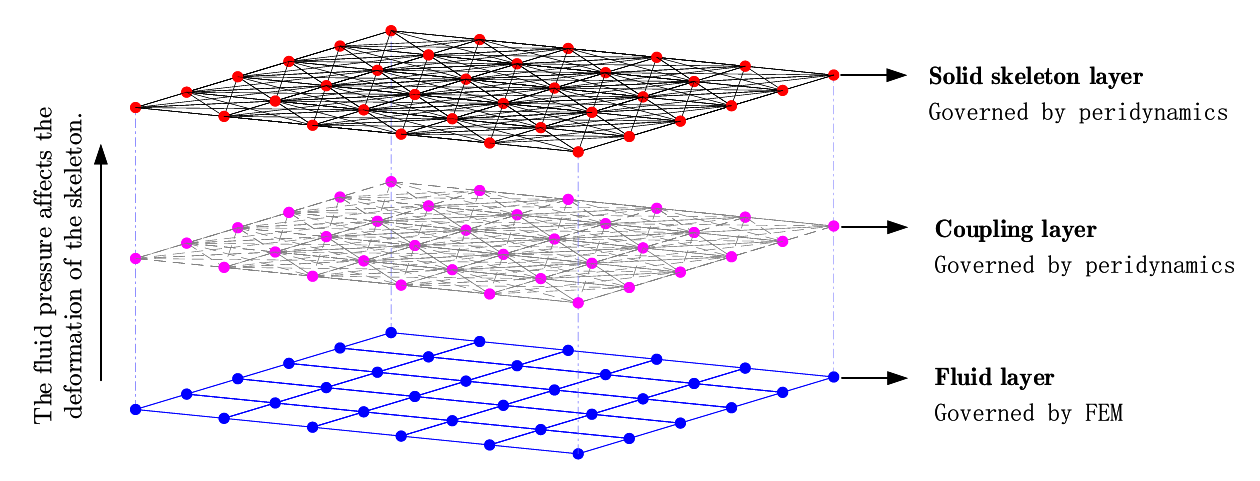}	\label{fig3_2:sub1}}
\\
\subfloat[Effect of solid on
	fluid.]{\includegraphics[scale=0.55]{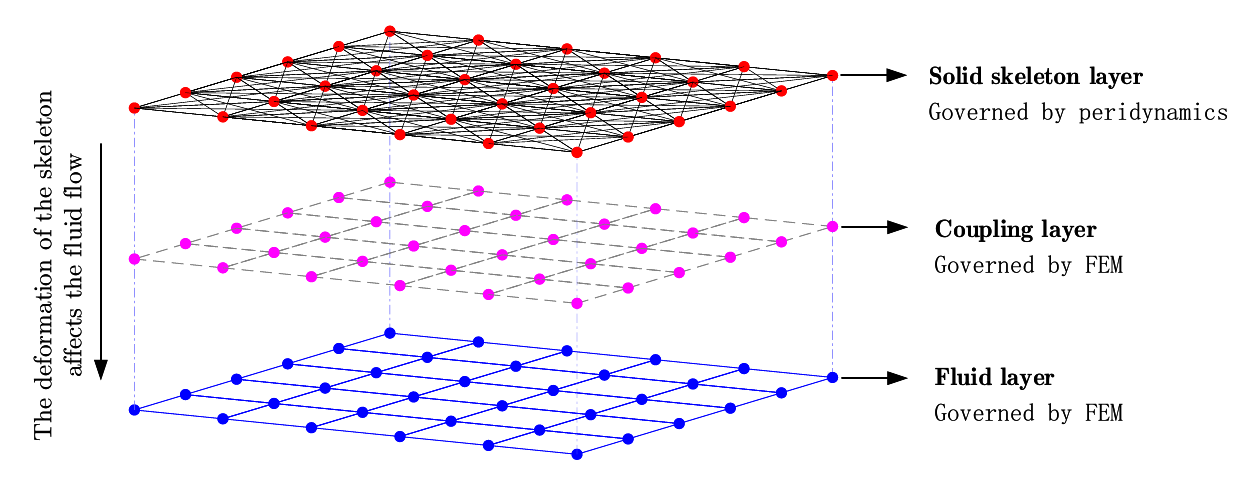}\label{fig3_2:sub2}}
\caption{Schematic diagram of bidirectional influence in the coupling
process, redrawn with permission from Ni et al. \citep{ni2020hybrid}.}
\label{coupling_scheme}
\end{figure}

Following the description of the hybrid modeling approach in %
\citep{ni2020hybrid}, the Galerkin finite element method is adopted for the
governing equation of the fluid flow, while the peridynamics for the
deformation of the solid phase. The fluid domain is discretized by using
4-node FE elements and the solid domain by using a regular PD grid. PD and
FE nodes share the same coordinates, but this is not required by our
approach. The implementation of the hydro-mechanical coupling process is
described in Fig. \ref{coupling_scheme}. 

\subsection{Discretization of governing equations}

The discretized peridynamic equation of motion of the current node $%
\boldsymbol{x}_{i}$ is written as:%
\begin{equation}
\begin{array}{l}
\rho \boldsymbol{\ddot{u}}_{i}^{t}=\sum_{j=1}^{N_{H_{i}}}\left\{ \text{ 
\textbf{\b{T}}}\left[ \boldsymbol{x}_{i},t\right] \left\langle \boldsymbol{%
\xi }_{ij}\right\rangle -\text{ \textbf{\b{T}}}\left[ \boldsymbol{x}_{j},t%
\right] \left\langle \boldsymbol{-\xi }_{ij}\right\rangle \right\} \cdot
V_{j} \\ 
-2\alpha \sum_{j=1}^{N_{H_{i}}}\left[ p_{i}\frac{\underline{\mathit{w}}%
\left\langle \boldsymbol{\xi }_{ij}\right\rangle \text{ }\underline{x}%
\left\langle \boldsymbol{\xi }_{ij}\right\rangle }{m\left( \boldsymbol{x}%
_{i}\right) }\underline{\boldsymbol{M}}\left\langle \boldsymbol{\xi }%
_{ij}\right\rangle -\text{ }p_{j}\frac{\underline{\mathit{w}}\left\langle 
\boldsymbol{\xi }_{ij}\right\rangle \text{ }\underline{x}\left\langle 
\boldsymbol{\xi }_{ij}\right\rangle }{m\left( \boldsymbol{x}_{j}\right) }%
\underline{\boldsymbol{M}}\left\langle -\boldsymbol{\xi }_{ij}\right\rangle %
\right] \cdot V_{j}+\boldsymbol{b}_{i}^{t}%
\end{array}
\label{3.1.10}
\end{equation}%
where $N_{H_{i}}$ is the number of family nodes of $\boldsymbol{x}_{i}$, $%
\boldsymbol{x}_{j}$\ represents $\boldsymbol{x}_{i}$'s family node, $V_{j}$\ is the
volume of node $\boldsymbol{x}_{j}$.

Under the assumption of small deformation, the global form of Eq.(\ref%
{3.1.10}) can be written as:%
\begin{equation}
\boldsymbol{M}^{PD}\boldsymbol{\ddot{u}}+\boldsymbol{K}^{PD}\boldsymbol{u}-%
\boldsymbol{Q}^{PD}p=f^{PD}  \label{3.1.11}
\end{equation}%
in which $\boldsymbol{M}^{PD}$, $\boldsymbol{K}^{PD}$ and $\boldsymbol{Q}%
^{PD}$\ are the mass, stiffness and \textquotedblleft
coupling\textquotedblright\ matrices of the PD domain, respectively. Note
that $\boldsymbol{M}^{PD}$\ is usually taken as a lumped mass matrix. The
method of obtaining $\boldsymbol{K}^{PD}$ of OSB-PD equations can be found
in \citep{sarego2016linearized}. Assuming that $\underline{\boldsymbol{M}}%
\left\langle \boldsymbol{\xi }_{ij}\right\rangle =\left[ M_{x},M_{y}\right] $%
, then the \textquotedblleft coupling\textquotedblright\ matrix for $%
\boldsymbol{\xi }_{ij}$ can be given as:%
\begin{equation}
\boldsymbol{Q_{ij}^{PD}=}2\alpha \underline{\mathit{w}}\left\langle 
\boldsymbol{\xi }_{ij}\right\rangle \underline{x}\left\langle \boldsymbol{\
\xi }_{ij}\right\rangle V_{i}V_{j}\left[ 
\begin{array}{cccc}
\frac{M_{x}}{m\left( \boldsymbol{x}_{i}\right) } & \frac{M_{y}}{m\left( 
\boldsymbol{x}_{i}\right) } & -\frac{M_{x}}{m\left( \boldsymbol{x}%
_{i}\right) } & -\frac{M_{y}}{m\left( \boldsymbol{x}_{i}\right) } \\ 
\frac{M_{x}}{m\left( \boldsymbol{x}_{j}\right) } & \frac{M_{y}}{m\left( 
\boldsymbol{x}_{j}\right) } & -\frac{M_{x}}{m\left( \boldsymbol{x}%
_{j}\right) } & -\frac{M_{y}}{m\left( \boldsymbol{x}_{j}\right) }%
\end{array}%
\right] ^{T}  \label{3.1.12}
\end{equation}%
where the influence function $\underline{\mathit{w}}\left\langle \boldsymbol{%
\ \xi }_{ij}\right\rangle $ used for the coupling bonds should be specified
as $\underline{\mathit{w}}=1$.

The discretized governing equations of the fluid flow assume the following
form \citep{lewis1998finite}:%
\begin{equation}
\boldsymbol{S\dot{p}}+\boldsymbol{Q}^{T}\boldsymbol{\dot{u}}+\boldsymbol{Hp=q%
}^{w}  \label{3.1.1}
\end{equation}%
where $\boldsymbol{Q}$ is the coupling matrix, $\boldsymbol{H}$ is the
permeability matrix, $\boldsymbol{S}$ is the compressibility matrix, all of
them can be obtained by assembling the corresponding element matrices. When
using the shape functions $\boldsymbol{N}_{u}$ (for displacement) and $%
\boldsymbol{N}_{p}$ (for pressure), the matrices in Eq.(\ref{3.1.1}) can be
given as \citep{milanese2016avalanches,peruzzo2019dynamics}:%
\begin{equation}
\boldsymbol{Q=}\int\nolimits_{\Omega }\left( \boldsymbol{LN}_{u}\right)
^{T}\alpha \boldsymbol{mN}_{p}\text{d}\Omega  \label{3.1.2}
\end{equation}%
\begin{equation}
\boldsymbol{H=}\int\nolimits_{\Omega }\left( \nabla \boldsymbol{N}%
_{p}\right) ^{T}\frac{k}{\mu _{w}}\left( \nabla \boldsymbol{N}_{p}\right) 
\text{d}\Omega  \label{3.1.3}
\end{equation}%
\begin{equation}
\boldsymbol{S=}\int\nolimits_{\Omega }\boldsymbol{N}_{p}^{T}s\boldsymbol{N}%
_{p}\text{d}\Omega  \label{3.1.4}
\end{equation}%
in which $\boldsymbol{L}$ is the differential operator defined as: 
\begin{equation}
\boldsymbol{L=}\left[ 
\begin{array}{ccc}
\frac{\partial }{\partial x} & 0 & \frac{\partial }{\partial y} \\ 
0 & \frac{\partial }{\partial y} & \frac{\partial }{\partial x}%
\end{array}%
\right] ^{T}  \label{3.1.5}
\end{equation}%
$\boldsymbol{m}$\ is a vector given as: 
\begin{equation}
\boldsymbol{m}=\left[ 1,1,0\right] ^{T}  \label{3.1.9}
\end{equation}

As described in $sect.$ $3.3$, the volumetric coupling term $\left( \rho
\alpha \frac{\partial \varepsilon _{v}}{\partial t}\right) $ will be removed
from the governing equation for the flow in the fracture domain.

\subsection{Staggered solution scheme of the hydro-mechanical coupled system}

As in \citep{lewis1998finite,zienkiewicz2000finite}, Eqs. \ref{3.1.11} and %
\ref{3.1.1} can be combined and rewritten in the following format: 
\begin{equation}
\left[ 
\begin{array}{cc}
\boldsymbol{M}^{PD} & \boldsymbol{0} \\ 
\boldsymbol{0} & \boldsymbol{0}%
\end{array}%
\right] \left[ 
\begin{array}{c}
\boldsymbol{\ddot{u}} \\ 
\boldsymbol{\ddot{p}}%
\end{array}%
\right] +\left[ 
\begin{array}{cc}
\boldsymbol{0} & \boldsymbol{0} \\ 
\boldsymbol{Q}^{T} & \boldsymbol{S}%
\end{array}%
\right] \left[ 
\begin{array}{c}
\boldsymbol{\dot{u}} \\ 
\boldsymbol{\dot{p}}%
\end{array}%
\right] +\left[ 
\begin{array}{cc}
\boldsymbol{K}^{PD} & \boldsymbol{-Q}^{PD} \\ 
\boldsymbol{0} & \boldsymbol{H}%
\end{array}%
\right] \left[ 
\begin{array}{c}
\boldsymbol{u} \\ 
\boldsymbol{p}%
\end{array}%
\right] =\left[ 
\begin{array}{c}
\boldsymbol{f} \\ 
\boldsymbol{q}^{w}%
\end{array}%
\right]  \label{coupled_matrix_eq}
\end{equation}

\textquotedblleft Monolithic\textquotedblright\ and \textquotedblleft
staggered\textquotedblright\ algorithms are commonly adopted to solve such a
coupled system. Details of the \textquotedblleft
monolithic\textquotedblright\ solution algorithm can be found in %
\citep{zienkiewicz2000finite,peruzzo2019dynamics,ni2020hybrid}. It is
usually both memory- and time-intensive when used for the solution of hydro-mechanical coupled
problems involving fracture advancement. Therefore often, a \textquotedblleft
staggered\textquotedblright\ algorithm is preferred to obtain the dynamic
solutions of such kind of coupled problems.

In a \textquotedblleft staggered\textquotedblright\ algorithm, the first row
and second row in Eq. \ref{coupled_matrix_eq} are solved sequentially, and
the previously solved field variables ($\boldsymbol{u}$ and $\boldsymbol{p}$%
) are used as conditions of the next solving sequence. For the described
hydro-mechanical coupled problems, the two steps in each solving sequence
can be summarized as follows:

--step 1: solve the pressure field ($\boldsymbol{p}^{n+1}$) using the
following implicit time integration iteration \citep{ni2020hybrid}:%
\begin{equation}
\boldsymbol{p}^{n+1}=\left[ \boldsymbol{S}+\vartheta \Delta t\boldsymbol{H}%
\right] ^{-1}\left\{ \left[ \boldsymbol{S}-\left( 1-\vartheta \right) \Delta
t\boldsymbol{H}\right] \boldsymbol{p}^{n}-\Delta t\boldsymbol{q}^{w}+%
\boldsymbol{Q}^{T}\left( \mathbf{u}^{n}-\mathbf{u}^{n-1}\right) \right\}
\label{3.4.5}
\end{equation}

--step 2: solve the displacement field ($\boldsymbol{u}^{n+1}$) of the solid
domain by using appropriate time integration algorithm.

This paper aims to capture the dynamic phenomena in hydraulic fracture
problems. Thus, a modified explicit central difference time integration
scheme proposed in \citep{taylor1989pronto} is adopted here to obtain the
dynamic solution of the OSB-PD model. In the chosen time integration scheme,
the velocities are integrated with a forward difference and the
displacements with a backward difference. The velocity and displacement at
the $\left( n+1\right) ^{th}$ time increment can be obtained as:%
\begin{equation}
\begin{array}{l}
\boldsymbol{\dot{u}}^{n+1}=\boldsymbol{\dot{u}}^{n}+\Delta t\boldsymbol{%
\ddot{u}}^{n} \\ 
\boldsymbol{u}^{n+1}=\boldsymbol{u}^{n}+\Delta t\boldsymbol{\dot{u}}^{n+1}%
\end{array}
\label{3.11}
\end{equation}%
where $\boldsymbol{\ddot{u}}^{n}$ is the acceleration at $n^{th}$\ time
increment and can be determined by using Newton's second law:%
\begin{equation}
\boldsymbol{\ddot{u}}^{n+1}=\mathbf{M}^{-1}\left( \mathbf{F}^{ext}-\mathbf{F}%
^{int}-\boldsymbol{Q}^{PD}\boldsymbol{p}\right)  \label{3.10}
\end{equation}%
where $\mathbf{F}^{ext}$ and $\mathbf{F}^{int}$ are the external and
internal force vectors, respectively, $\mathbf{M}$ is the diagonal mass
matrix. The $\Delta t$ in the above equations is the constant time
increment, and the explicit method for the undamped system requires a
critical time step for numerical stability. According to %
\citep{zhou2016numerical}, the stable time increment for PD model can be
defined as:%
\begin{equation}
\Delta t<\delta /c^{\prime }  \label{3.13}
\end{equation}%
where $c^{\prime }=\sqrt{\left( \lambda +2\mu \right) /\rho }$ is the
dilatational wave speed and $\lambda $ and $\mu $ are the Lame's elastic
constants of the material.

Note that the permeability and storage matrices ($\boldsymbol{H}$ and $%
\boldsymbol{S}$) need to be updated accordingly in each time step when
involving cracks.

\textcolor{blue}{\textbf{Remark 1.} The FE formulation for fluid flow, even if local, contains de facto a length scale due to the fact that we have the natural presence of a gradient term through Darcy's Law and this inside a divergence operator. The result is a Laplace's operator which is known to introduce a length scale even if sometimes weak. This fact has been extensively investigated in dynamic strain localization analysis \citep{zhang1999interal}. As stated in \citep{ni2020hybrid}, when solving the hybrid model, the coupling matrices in the first row and second row of Eq. (\ref{coupled_matrix_eq}) should be generated non-locally and locally, respectively. However, because of the inconsistency of the length scales in the hybrid local and non-local coupled system, it has been found that it may entail problems in the simulation of dynamic problems. It is usually difficult to ensure the consistency of length scales in the both domains. Thus, although not theoretically rigorous, it is a convenient choice to adopt the same formulation of coupling matrices to avoid numerical instability. In addition, recently some criticism has been put forward about wave dispersion characteristics of PD models in \citep{bavzant2016wave}, which could be another factor affecting the dynamic solutions of the coupled system and deserves further scrutiny.}

\section{Wave propagation in one-dimensional dynamic consolidation problem:
a benchmark case}

\begin{figure}[h]
\begin{center}
\includegraphics[scale=0.55]{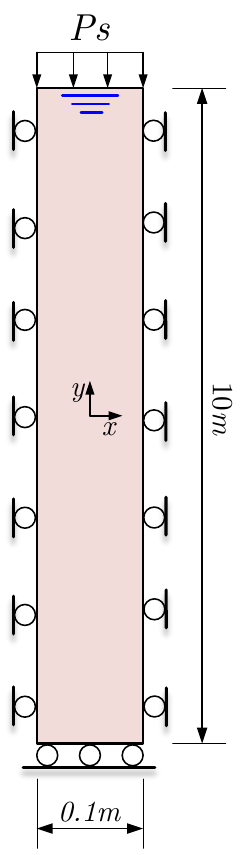}
\end{center}
\caption{Geometry and boundary conditions of the one-dimensional dynamic
consolidation problem.}
\label{fig4_1}
\end{figure}

In this section, an example of one-dimensional dynamic consolidation is
presented to validate the proposed approach. This example was solved by
Schanz and Cheng \citep{schanz2000transient} for the transient wave
propagation of the displacement and pore pressure during the dynamic
consolidation process and analytical solutions were obtained by using the
Dubner and Abate's method \citep{dubner1968numerical} for comparison. In
order to reproduce this example and compare our solution with the analytical
one, the behavior of a poroelastic column shown in Fig. \ref{fig4_1} is
simulated under the following boundary conditions. The top is drained with a
pore pressure of $p=0$ and subjected to a sudden surface pressure $P_{s}=1Pa$%
, while the other boundaries are impermeable and constrained in their normal
direction. The mechanical and fluid parameters used in the simulations are
listed in Tab. \ref{Tab.1}. Referring to \citep{schanz2000transient}, a
fixed time step of $\Delta t=1\times 10^{-5}s$ and $\vartheta =0.5$ are
chosen for the time integration. 
\begin{table}[h]
\caption{Mechanical and fluid parameters used in the example of dynamic
consolidation problem.}
\label{Tab.1}\centering
\begin{tabular}{cccccc}
\toprule $E$ & $\nu$ & $\rho_{r}$ & $\rho_{w}$ & $\alpha$ &  \\ \hline
$0.254 GPa$ & $0.3$ & $2700kg/m^{3}$ & $1000kg/m^{3}$ & $0.7883$ &  \\ \hline
$n$ & $K_{r}$ & $K_{w}$ & $\mu$ & $k$ &  \\ \hline
$0.48$ & $11 GPa $ & $3.3 GPa$ & $1\times10^{-3}Pa\cdot s$ & $%
3.55\times10^{-12}m^{2}$ &  \\ 
\bottomrule &  &  &  &  & 
\end{tabular}
%	}
\end{table}

\begin{figure}[h!]
\centering  
\subfloat[Variation of vertical displacement on the top edge versus
time.]{\includegraphics[scale=0.65]{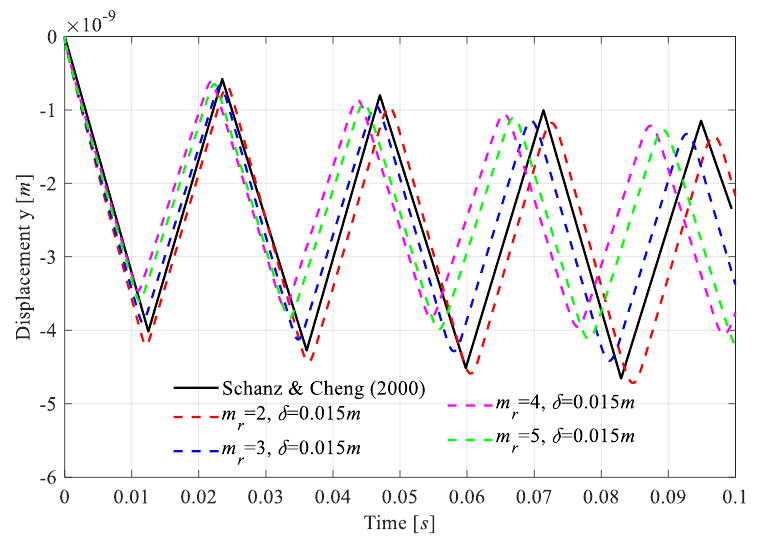}\label{fig5_1:sub1}}%
\subfloat[Variation of pore pressure on the bottom edge versus
time.]{\includegraphics[scale=0.65]{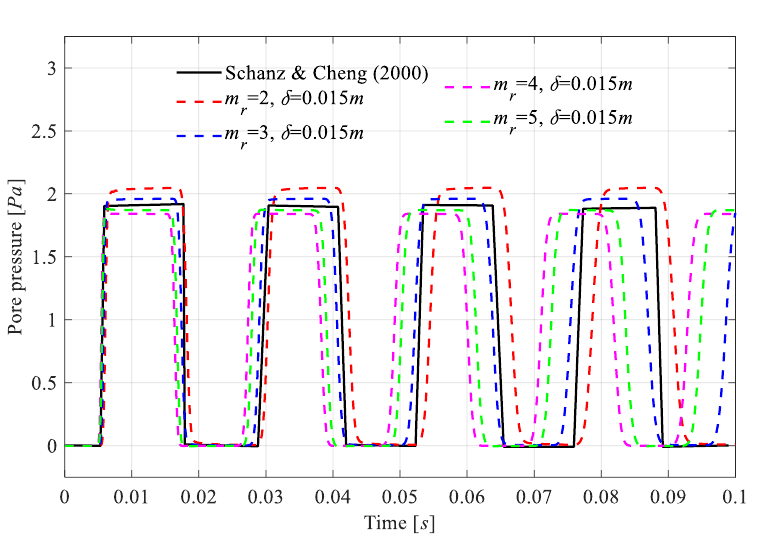}	\label{fig5_1:sub2}}
\caption{$m_{r}$-convergence study of the dynamic consolidation example with $%
\protect\delta=0.015m$.}
\label{fig5_1}
\end{figure}

\begin{figure}[h!]
\centering  
\subfloat[Variation of vertical displacement on the top edge versus
time.]{\includegraphics[scale=0.65]{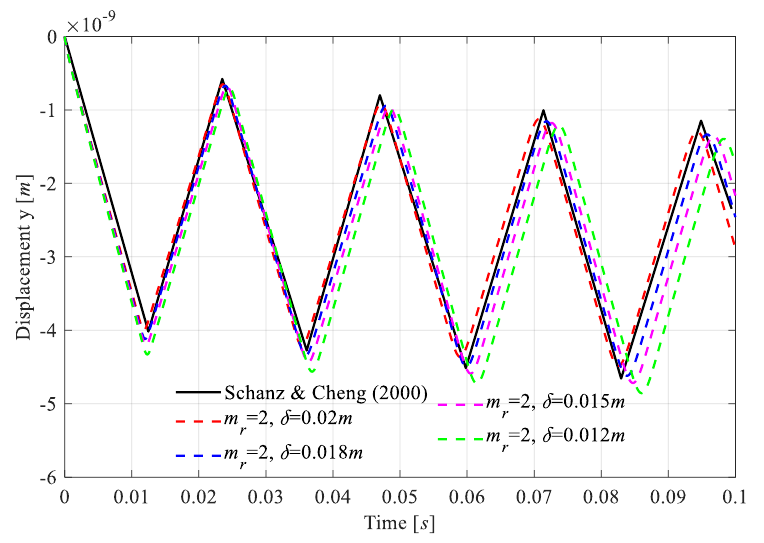}\label{fig5_2:sub1}}
\subfloat[Variation of pore pressure on the bottom edge versus
time.]{\includegraphics[scale=0.65]{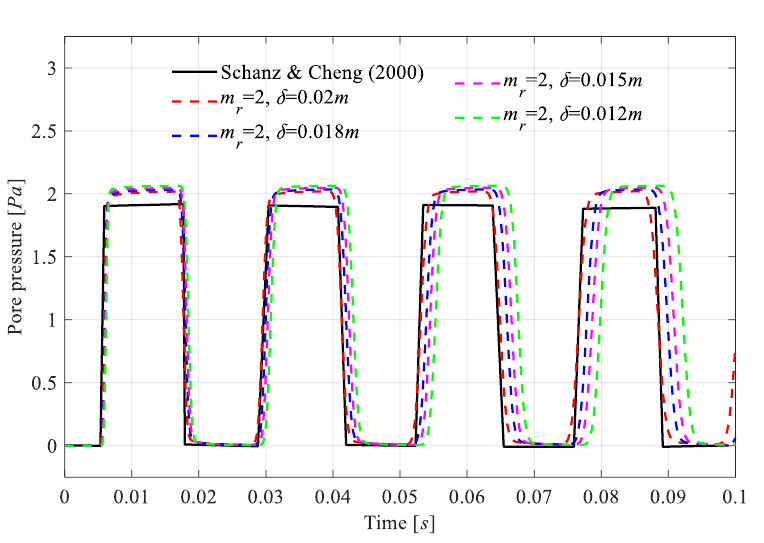}	\label{fig5_2:sub2}}
\caption{$\protect\delta$-convergence study of the dynamic consolidation
example with $m_{r}=2$.}
\label{fig5_2}
\end{figure}

Discretization parameters are the key factors affecting the solutions of
peridynamic models \citep{le2014two}. It is therefore necessary to perform
convergence studies on discretization parameters before using
peridynamic-based tools for numerical analysis. $\delta $-convergence and $%
m_{r}$-convergence studies are hence performed here to identify the most
suitable discretization parameters in terms of accuracy and computational
cost. In the $m_{r}$-convergence study, the $\delta =0.015m$ and $%
m_{r}=2,3,4,5$ are adopted respectively. The variations of the vertical
displacement on top edge and pore pressure on the bottom edge are plotted in
Figs. \ref{fig5_1:sub1} and \ref{fig5_1:sub2}, respectively. The comparison
between the numerical and analytical solutions suggests that $m_{r}=2$ is an
acceptable choice for such a hydro-mechanical coupled problem if both the
accuracy of displacement and pore pressure are considered. Consequently, $%
m_{r}=2$ in conjunction with $\delta =0.012m,0.015m,0.018m,0.02m$ are chosen
for the $\delta $-convergence study. The numerical solutions of the $\delta $%
-convergence study are plotted in Figs. \ref{fig5_2:sub1} and \ref%
{fig5_2:sub2}. According to the numerical results, we can conclude that the
proposed hybrid modeling approach can properly simulate the dynamic response
in saturated porous media. However, the discretization parameters of
peridynamic models have a remarkable influence on simulation results. The proposed PD-based hybrid model shows a different convergence performance in such a hydro-mechanical coupled problem from that in mechanical-only problems \citep{le2014two}. In general, $%
m_{r}=2$ is an economical choice to obtain acceptable results with the
hybrid model when choosing the appropriate grid size.

\newpage

\section{Dynamic forerunning fracture in porous structures}

\begin{figure}[h!]
\centering  
\subfloat[Geometry and
constraints.]{\includegraphics[scale=0.55]{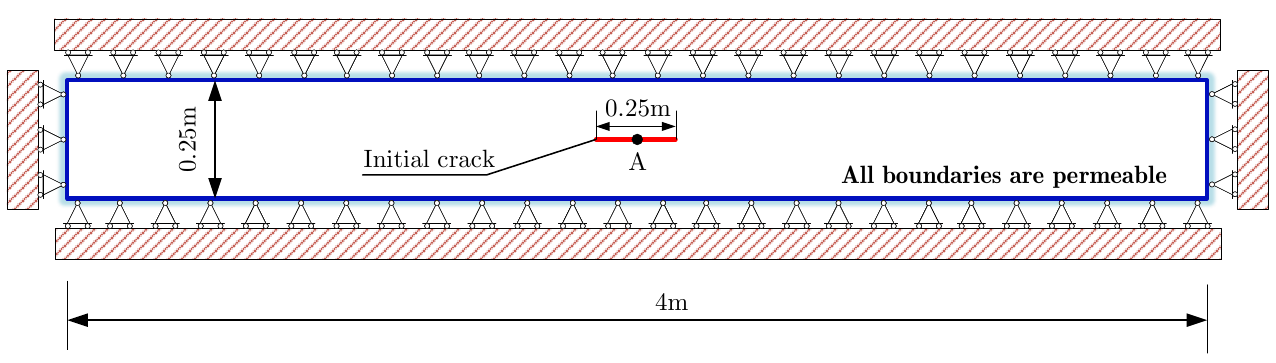}\label{fig6_1_1:sub1}}%
\\
\subfloat[Loading
conditions.]{\includegraphics[scale=0.45]{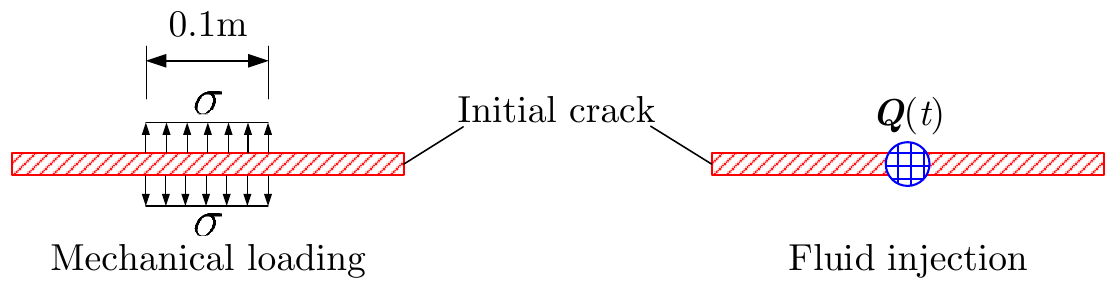}\label{fig6_1_1:sub2}}
\caption{Geometry and boundary conditions of the rectangular structure with
a central initial crack.}
\label{fig6_1_1}
\end{figure}

In \citep{cao2018porous,peruzzo2019dynamics,peruzzo2019stepwise}, a
cantilever beam sample was adopted to study the stepwise crack advancement
and fluid pressure oscillations in porous media fracturing dynamics. In %
\citep{peruzzo2019stepwise} as well as in \citep{slepyan2015forerunning},
fracture advancement with forerunning was observed in the dry case. Inspired
by that, a rectangular structure is solved in this section by using the
proposed approach for investigating dynamic fracturing with possible
forerunning in both dry and fully saturated porous media. The geometry and
constraints of the selected structure are drawn in Fig. \ref{fig6_1_1:sub1}.
Two different loading conditions, mechanical loading and fluid injection,
are adopted as respectively shown in Fig. \ref{fig6_1_1:sub2}. The
mechanical and fluid parameters used in the simulations are listed in Tab. %
\ref{Tab.2}. Note that the fracture energy release rate is taken as $%
G_{c}=1J/m^{2}$ for making it easier to achieve fracture advancement in the
structure. Since the described problem is symmetric, only the right half of the structure is modeled in the simulation. The whole domain is discretized by quadrilateral meshes with
grid size of $\Delta x=2.5\times 10^{-3}m$. $m_{r}=2$ is used for the
peridynamic model and thus the horizon radius is $\delta=2\times\Delta
x=5\times 10^{-3}m$. A fixed time step of $\Delta t=1\times10^{-7}s$ and 
$\vartheta=0.5$ are used for the time integration. The simulation time
duration is chosen as $1\times10^{-3}s$. \textcolor{blue}{The boundary conditions of the solid and fluid portions are applied non-locally and locally, respectively. The method of applying boundary conditions to peridynamic model can be found in \citep{ni2018peridynamic}, and to FEM model in \citep{lewis1998finite}.} According to the two loading
conditions and the material saturation, three different cases are solved. \textcolor{blue}{For the purpose of guiding the wave propagation, the crack is only allowed to propagate along the horizontal centerline of the specimen.}

\begin{table}[h!]
\caption{Mechanical and fluid parameters used in the forerunning fracturing
examples.}
\label{Tab.2}\centering
\begin{tabular}{cccccc}
\toprule $E$ & $\nu$ & $G_{c}$ & $\rho_{r}$ & $\rho_{w}$ &  \\ \hline
$10 GPa$ & $0.25$ & $1J/m^{2}$ & $1000kg/m^{3}$ & $1000kg/m^{3}$ &  \\ \hline
$\alpha$ & $n$ & $K_{w}$ & $\mu$ & $k$ &  \\ \hline
$1$ & $0.002 $ & $2.2 GPa$ & $1\times10^{-3}Pa\cdot s$ & $1
\times10^{-12}m^{2}$ &  \\ 
\bottomrule &  &  &  &  & 
\end{tabular}
%	}
\end{table}

\begin{figure}[h!]
\centering  
\subfloat[Damage levels at $2.77\times
10^{-4}s$.]{\includegraphics[scale=0.7]{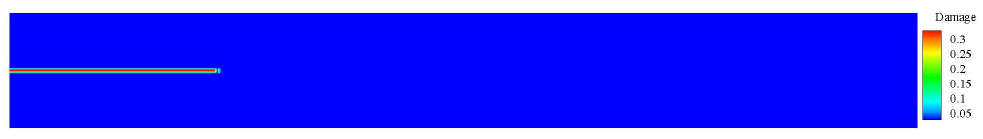}\label{fig6_1:sub1}}%
\\
\subfloat[Damage levels at $4.42\times
10^{-4}s$.]{\includegraphics[scale=0.7]{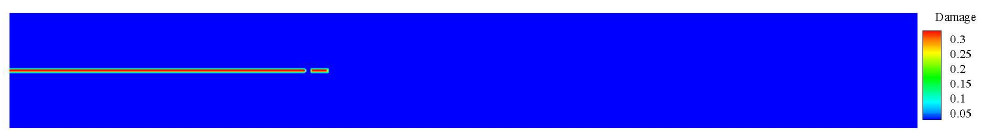}	\label{fig6_1:sub2}}%
\\
\subfloat[Damage levels at $8\times
10^{-4}s$.]{\includegraphics[scale=0.7]{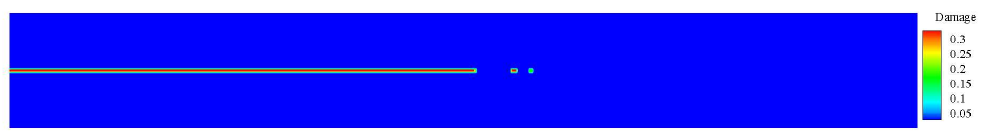}\label{fig6_1:sub3}}%
\\
\subfloat[Damage levels at $1\times
10^{-3}s$.]{\includegraphics[scale=0.7]{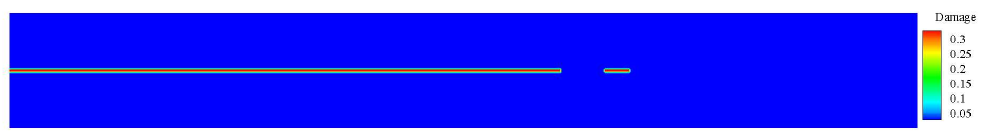}	\label{fig6_1:sub4}}
\caption{Forerunning fracture in the front of the crack tip in dry structure
under mechanical loading condition.}
\label{fig6_1}
\end{figure}

The first case (case 1) considers a dry porous material under mechanical
loading condition. As shown in Fig. \ref{fig6_1_1:sub2}, a constant pressure 
$\boldsymbol{\sigma}=15 MPa$ is applied on the surface of the initial crack
near the central position to force crack opening. Under the action of the
external force, the crack propagates from the initial crack tip. As shown in
Figs. \ref{fig6_1:sub1} to \ref{fig6_1:sub4}, forerunning events are
observed during the simulation. 

\begin{figure}[h!]
\centering  
\subfloat[Damage levels at $4\times
10^{-5}s$.]{\includegraphics[scale=0.7]{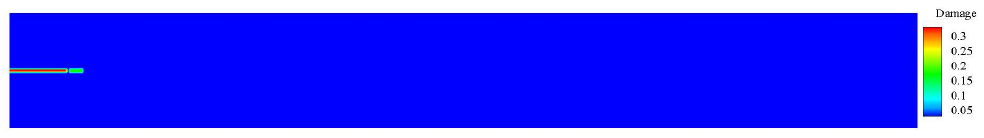}\label{fig6_2:sub1}}%
\\
\subfloat[Damage levels at $7.56\times
10^{-4}s$.]{\includegraphics[scale=0.7]{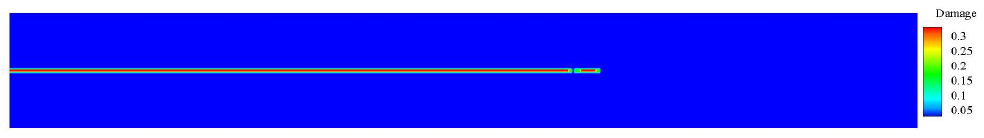}	\label{fig6_2:sub2}}%
\\
\subfloat[Damage levels at $9.41\times
10^{-4}s$.]{\includegraphics[scale=0.7]{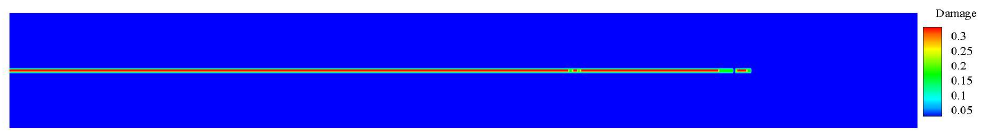}\label{fig6_2:sub3}}%
\\
\subfloat[Damage levels at $1\times
10^{-3}s$.]{\includegraphics[scale=0.7]{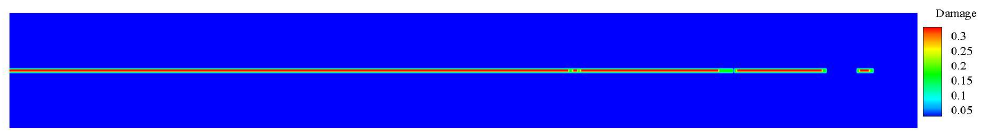}	\label{fig6_2:sub4}}
\caption{Forerunning fracture in the front of the crack tip in saturated
structure under mechanical loading condition.}
\label{fig6_2}
\end{figure}

\begin{figure}[h!]
\centering  
\subfloat[Distribution of pore pressure at $2.5\times
10^{-4}s$.]{\includegraphics[scale=0.7]{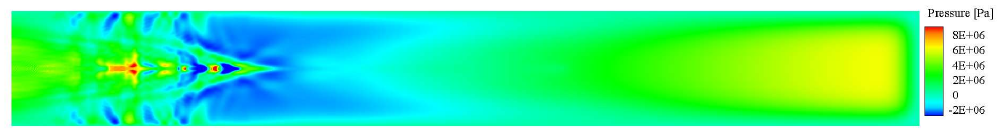}\label{fig6_3:sub1}}%
\\
\subfloat[Distribution of pore pressure at $2.6\times
10^{-4}s$.]{\includegraphics[scale=0.7]{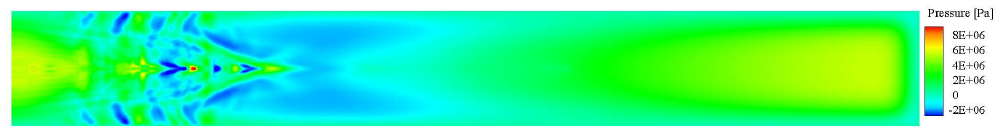}	\label{fig6_3:sub2}}%
\\
\subfloat[Distribution of pore pressure at $2.7\times
10^{-4}s$.]{\includegraphics[scale=0.7]{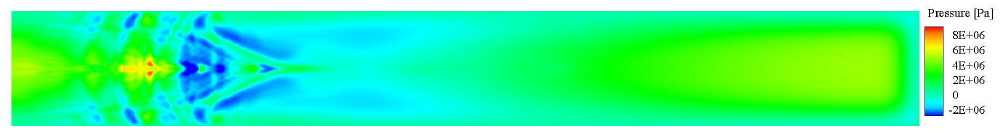}\label{fig6_3:sub3}}%
\\
\subfloat[Distribution of pore pressure at $2.8\times
10^{-4}s$.]{\includegraphics[scale=0.7]{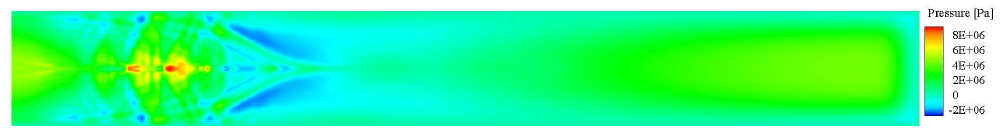}	\label{fig6_3:sub4}}
\caption{Wave propagation of pore pressure in saturated structure under
mechanical loading condition.}
\label{fig6_3}
\end{figure}

The second case (case 2) is carried out with a saturated porous material
under mechanical loading condition by applying the same constant pressure on
the surface of the initial crack. Under the action of applied external
force, the initial crack opens gradually, the material domain deforms and
the pore pressure in the whole domain changes accordingly. In this case, the
crack propagation occurs under the combined action of external force and
pore pressure. As shown in Figs. \ref{fig6_2:sub1} to \ref{fig6_2:sub4},
forerunning events also exist in this case. The distributions of pore
pressure at $2.5\times10^{-4}s$, $2.6\times10^{-4}s$, $2.7\times10^{-4}s$
and $2.8\times10^{-4}s$ are plotted in Figs. \ref{fig6_3:sub1} to \ref%
{fig6_3:sub4}, respectively. Changes shown in these images reveal the wave
propagation of pore pressure in the simulated structure. It is worth to
observe by comparing the fracture position reached after $1\times10^{-3}s$
that the fracture propagates faster in the saturated medium thanks to the
combined action of stress and pressure waves, and this result can also be
inferred from Fig. \ref{fig6_6}. 

\begin{figure}[h!]
\centering  
\subfloat[Damage levels at $1.54\times
10^{-4}s$.]{\includegraphics[scale=0.7]{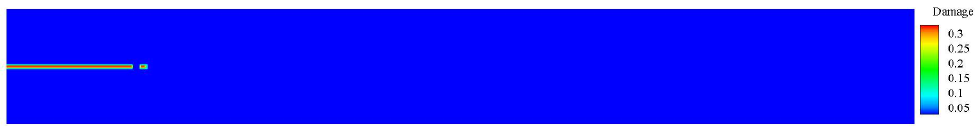}\label{fig6_4:sub1}}%
\\
\subfloat[Damage levels at $3.38\times
10^{-4}s$.]{\includegraphics[scale=0.7]{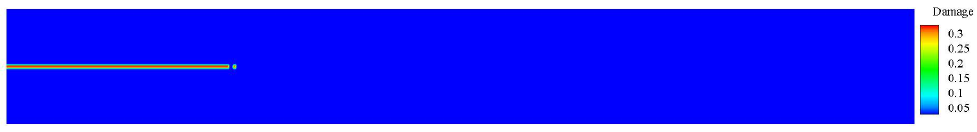}	\label{fig6_4:sub2}}%
\\
\subfloat[Damage levels at $5.51\times
10^{-4}s$.]{\includegraphics[scale=0.7]{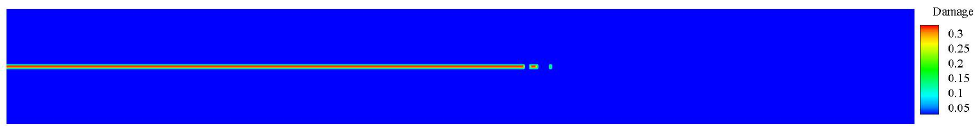}\label{fig6_4:sub3}}%
\\
\subfloat[Damage levels at $7.23\times
10^{-4}s$.]{\includegraphics[scale=0.7]{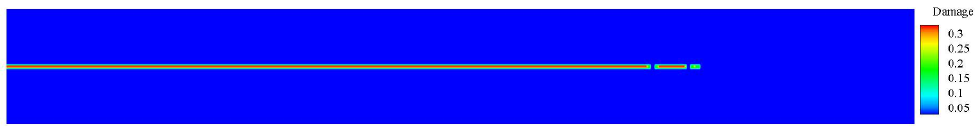}	\label{fig6_4:sub4}}
\caption{Forerunning fracture in the front of the crack tip in saturated
structure under fluid injection.}
\label{fig6_4}
\end{figure}

\begin{figure}[h!]
\centering  
\subfloat[Distribution of pore pressure at $2.3\times
10^{-4}s$.]{\includegraphics[scale=0.7]{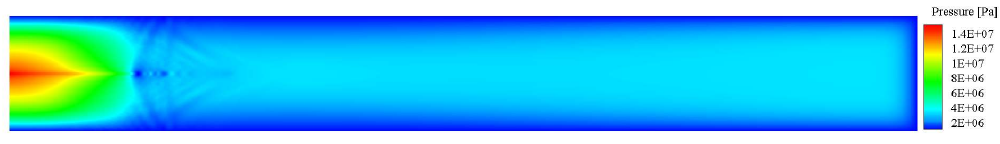}\label{fig6_5:sub1}}%
\\
\subfloat[Distribution of pore pressure at $2.4\times
10^{-4}s$.]{\includegraphics[scale=0.7]{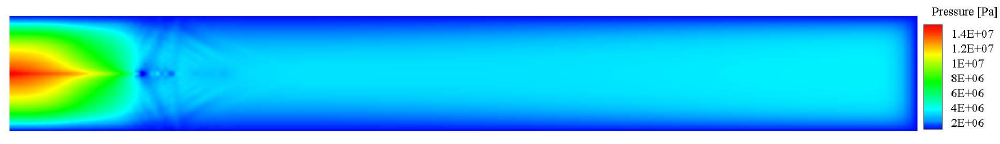}	\label{fig6_5:sub2}}%
\\
\subfloat[Distribution of pore pressure at $2.5\times
10^{-4}s$.]{\includegraphics[scale=0.7]{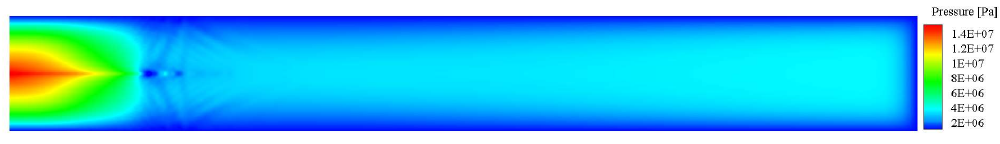}\label{fig6_5:sub3}}%
\\
\subfloat[Distribution of pore pressure at $2.6\times
10^{-4}s$.]{\includegraphics[scale=0.7]{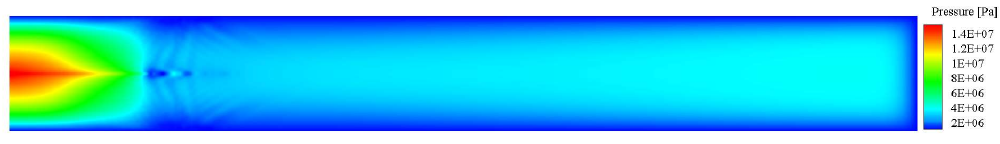}	\label{fig6_5:sub4}}
\caption{Wave propagation of pore pressure in saturated structure under
fluid injection.}
\label{fig6_5}
\end{figure}

\begin{figure}[h!]
\begin{center}
\includegraphics[scale=0.8]{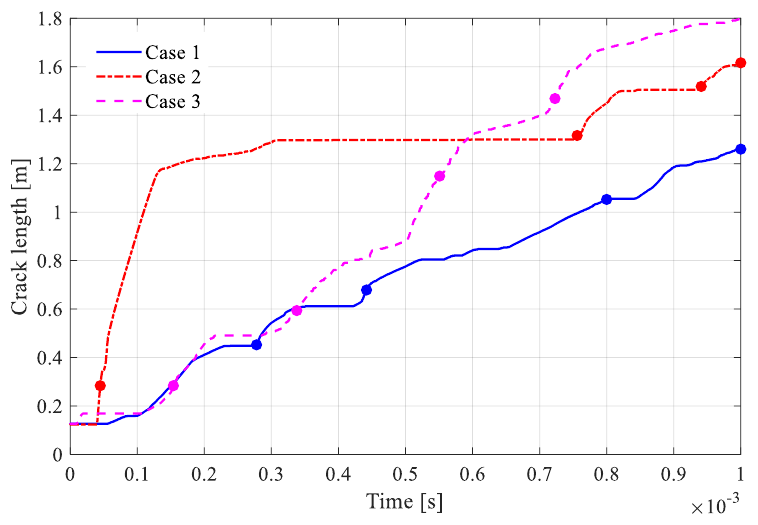}
\end{center}
\caption{Variation of crack length with time; case 1 is a dry structure under mechanical loading, case 2 a saturated structure under mechanical loading and case 3 a saturated structure with fluid injection.
}
\label{fig6_6}
\end{figure}

The last case (case 3) is carried out with the saturated porous material
under the fluid injection. The fluid is injected at the center of initial
crack with a constant volume rate of $Q=1 m^{3}/s$. Figs. \ref{fig6_4:sub1}
to \ref{fig6_4:sub4} show the crack patterns in this case at four different
time instants during the fracture process, and cracks are also observed in
front of the tip of the continuous crack. The distributions of pore pressure
at four time instants ($2.3\times10^{-4}s$, $2.4\times10^{-4}s$, $%
2.5\times10^{-4}s$ and $2.6\times10^{-4}s$) are plotted in Figs. \ref%
{fig6_5:sub1} to \ref{fig6_5:sub4}, respectively. A dynamic phenomenon of
wave propagation of pore pressure is also observed. 

The variation of crack length with time in the three cases is plotted in
Fig. \ref{fig6_6}. The solid dots marked in Fig. \ref{fig6_6} correspond to
the time instants in Figs. \ref{fig6_1}, \ref{fig6_2} and \ref{fig6_4},
respectively. Based on the locations of these solid dots and their corresponding distributions of damage levels, it can be found
that there appears a jump in crack advancement when forerunning fracture
occurs. All the examples shown in this section suggest numerically that the
stepwise advancement of fractures and the forerunning fracturing phenomenon
stressed in \citep{cao2018porous,peruzzo2019stepwise} do exist in the
dynamic fracture of dry and saturated porous media. 
\newpage
\section{Conclusions}
The hybrid FEM and Peridynamic modeling approach proposed in %
\citep{ni2020hybrid} is applied here to solve the problem of dynamic
fracture propagation in saturated porous media. In the hybrid modeling
approach, the fluid flow is governed by FE equations, while peridynamics is
used to solve the deformation of solids and capture the dynamic fracture
advancement. A staggered solution scheme is adopted to obtain the dynamic
solution of the hydro-mechanical coupled system. In each solving sequence,
an implicit time integration iteration from \citep{zienkiewicz2000finite} is
used to solve the FE equations of fluid flow, and a modified explicit
central difference time integration algorithm proposed in %
\citep{taylor1989pronto} is adopted for the peridynamic equations.

Firstly, an example of a one-dimensional dynamic consolidation problem is
solved by the presented approach as a benchmark case. Both $\delta-$%
convergence and $m_{r}-$convergence studies are performed to investigate the
influences of discretization parameters on the dynamic solutions of such a
coupled problem. Based on the comparison of numerical results and analytical
solutions, $m_{r}=2$ is suggested for the PD portion to obtain acceptable
dynamic solution of the described hydro-mechanical coupled system.
Subsequently, a series of numerical simulations of the dynamic crack
propagation in dry and fully saturated porous structures are carried out.
Forerunning fracture events, a phenomenon that has been observed
experimentally in geomaterials \citep{sammonds1989acoustic}, occur in these
simulations. Besides the mechanical wave, the wave propagation of pore
pressure has also been numerically confirmed to exist in saturated porous
structures. The interaction between mechanical waves and pore pressure waves
makes dynamic fracturing in saturated media more complicated, influencing
also the forerunning episodes. \textcolor{blue}{It is worth to point out that in \citep{Milanese2020forerunning} it has been shown that in dry bodies (i) forerunning is a undeniable source of stepwise fracturing advancement; (ii) forerunning is a stability problem; and (iii) forerunning increases the overall tip advancement speed and is, in fact, a mechanism for a crack to move faster when a steady-state propagation is no longer supported by the body due to a high level of external forces. It appears clearly from Fig. \ref{fig6_6} that the interaction between
stress and pressure waves in saturated media increases the propagation speed
of the fracture under mechanical loading even further when compared to dry specimens}, at
least when forerunning events happen. This observation may be important for
geophysics.

\clearpage

\section*{Acknowledgements}

T. Ni, U. Galvanetto and M. Zaccariotto acknowledge the support they received from
MIUR under the research project PRIN2017-DEVISU and from University of Padua
under the research project BIRD2018 NR.183703/18.

F. Pesavento would like to acknowledge the project 734370-BESTOFRAC
\textquotedblleft Environmentally best practices and optimisation in
hydraulic fracturing for shale gas/oil development\textquotedblright
-H2020-MSCA-RISE-2016 and the support he received from University of Padua
under the research project BIRD197110/19 \textquotedblleft Innovative models
for the simulation of fracturing phenomena in structural engineering and
geomechanics\textquotedblright .

B.A. Schrefler gratefully acknowledges the support of the Technische
Universität München - Institute for Advanced Study, funded by the German
Excellence Initiative and the TÜV SÜD Foundation. \clearpage

%\bibliographystyle{ametsoc2014}
%\bibliographystyle{unsrt} 
%\bibliography{mybib}

\end{document}